\documentclass[12pt]{amsart}
\usepackage{latexsym}
\usepackage{amsfonts}
\begin{document}
\thispagestyle{empty}

\newtheorem{prop}{Proposition}
\newtheorem{theo}{Theorem}
\newtheorem{lemm}{Lemma}

\title{The canonical solution operator to $\overline \partial $ restricted to 
Bergman spaces and spaces of entire functions}
\author{Friedrich Haslinger}

\address{Institut f\"ur Mathematik, Universit\"at Wien\\
Strudlhofgasse 4, A-1090 Wien, Austria\\
e-mail: friedrich.haslinger@univie.ac.at}

\begin{abstract}
In this paper we  obtain a necessary and
sufficient condition 
for the canonical solution operator to $\overline \partial $
restricted   to  radial symmetric Bergman spaces to be a Hilbert-Schmidt operator.
We also discuss compactness of the solution operator in spaces of entire functions
in one variable.
In the sequel we consider several examples and also treat the case
of weighted spaces of entire functions in several variables.
\end{abstract}

\maketitle

\section{Introduction}

In  this paper we consider Bergman spaces $A^2(D, d\mu),$ where
$D$  is  a  disc  in  $\mathbb C $ and $d\mu (z) = m(z)d\lambda
(z),$  $m$  being  a  radial  symmetric weight function. We also
suppose  that  the  monomials  $\{z^n  \},  \ n\in \mathbb N_0$
constitute an orthogonal basis in $A^2(D, d\mu).$

Let 
$$c_n^2= \int_D |z|^{2n}\,d\mu (z).$$

We solve the $\overline \partial $-equation $\overline \partial
u  =  g,$  where    $g \in  A^2(D, d\mu).$

It is pointed out in \cite{FS1} that in the proof that compactness of the solution
operator for $\overline \partial $ on $(0,1)$-forms implies 
that the boundary of $\Omega $
does not contain any analytic variety of dimension grearter than or equal to 1,
it is only used that there is a compact solution operator 
to $\overline \partial $ on the 
$(0,1)$-forms with holomorphic coefficients.
In this case compactness
of the solution operator restricted to $(0,1)$-forms with holomorphic
coefficients implies already compactness of the solution operator on
general $(0,1)$-forms.

A similar situation appears in \cite{SSU} where the Toeplitz $C^*$ -algebra
$\mathcal T (\Omega )$ is considered and the relation between the structure of 
$\mathcal T (\Omega )$ and the $\overline \partial $-Neumann problem is discussed 
(see \cite{SSU} , Corollary 4.6).

In many cases non-compactness of the canonical solution operator already happens
when the solution operator is restricted to the corresponding subspace of 
holomorphic functions (or $(0,1)$-forms with holomorphic coefficients, in the case
of several variables.)(see \cite{FS1} , \cite{SSU} , \cite{K2}). In this paper
we will show that this phenomenon also occurs in the Fock space in one variable.

The  question  of compactness of the $\overline \partial $-Neumann operator
is of interest for various reasons (see \cite{FS2} ).

We use the fact that the canonical solution operator $S$ to $\overline \partial $
restricted  to  $(0,1)$-forms with holomorphic coefficients can
be   expressed by an integral operator using the Bergman kernel
(see \cite{Has2}):
$$S(g)(z)=\int_D K(z,w)\, (\overline z -\overline w)\,
g(w)\,d\mu (w),$$

where  $g\in  A^2(D, d\mu )$ and $K(z,w)$ is the Bergman kernel
of $A^2(D, d\mu ).$

With  the  help  of  this    result  we  obtain a necessary and
sufficient condition 
for the canonical solution operator to $\overline \partial $
restricted   to  $A^2(D, d\mu)$ to be a Hilbert-Schmidt operator.
This  condition is expressed in terms of the sequence $(c_n)_n$
defined above.

It   turns   out  that  for  $D=  \mathbb  D$  and  $d\mu  (z)=
(1-|z|^2)^{\alpha  }  \,d\lambda  (z)  \  \alpha \ge 0$ 
the  canonical  solution operator is always Hilbert Schmidt, for
the case $\alpha =0 $ see \cite{Has2} .

In the following part we discuss compactness of the canonical solution operator
to $\overline \partial $ in spaces of entire functions in one complex variable. We show
that the canonical solution operator
for $\overline \partial $ as operator from $L^2(\mathbb C, e^{-|z|^2}) $ 
into itself is not compact.
This follows from the result that the canonical solution operator
for $\overline \partial $ restricted to weighted space of entire 
functions $A^2(\mathbb C, e^{-|z|^2}) $ 
(Fock space) into $L^2(\mathbb C, e^{-|z|^2}) $ already fails to be compact. 
Further it is
shown that the restriction to $A^2(\mathbb C, e^{-|z|^m}) \ ,\ m>2,$ is compact but not 
Hilbert Schmidt.

In  the  sequel we also consider the case of several complex
variables in a slightly different situation and show
that the canonical solution operator to
$\overline \partial  $  is  a  Hilbert-Schmidt operator for a wide class of
weighted spaces of entire functions using various methods from
abstract functional analysis (see \cite{MV}).

The canonical solution operator to $\overline \partial $
restricted   to  $(0,1)$-forms with holomorphic coefficients
can also be interpreted as the Hankel operator
$$H_{\overline z}(g)=(I-P)(\overline z g),$$
where  $P  :  L^2(\Omega )
\longrightarrow A^2(\Omega )$ denotes the Bergman projection.
See \cite{A}, \cite{AFP}, \cite{J},
\cite{W} and \cite{Z} for details.

\section{Radial symmetric Bergman spaces}

The canonical solution operator
has  the  properties    $\overline  \partial  S (g) = g$ and
$S(g) \perp A^2(D, d\mu ).$

\begin{theo}
The canonical solution operator
$$S : A^2(D, d\mu ) \longrightarrow L^2(D, d\mu )$$
is a Hilbert Schmidt operator if and only if
$$\lim_{n\to \infty}\frac{c_{n+1}^2}{c_n^2}< \infty ,$$
where
$$c_n^2= \int_D |z|^{2n}\,d\mu (z)$$.
\end{theo}

\begin{proof}
By \cite{MV} , 16.8,  we have to show that there exists a
complete  orthonormal  system  $(u_k  )_{k=0}^{\infty}$ 
of $A^2(D, d\mu )$ such that
$$\sum_{n=0}^{\infty }\|S(u_n )\|^2 < \infty .$$
For this purpose we take the complete orthonormal system
$u_k(z)=z^k/c_k  ,$  then  the  Bergman  kernel $K(z,w)$ can be
expressed in the form
$$K(z,w)=\sum_{k=0}^{\infty }\frac{z^k\overline w^k}{c_k^2}.$$
From \cite{Has2} we know that
$$S(u_n)(z)    =  \overline  z  u_n(z) - \int_D K(z,w)\overline
w \,u_n(w)\,d\mu (w).$$
First we calculate the integral
\begin{eqnarray*}
\int_D K(z,w)\overline w\, u_n(w)\,d\mu (w) & = &
\int_D \overline w\, \frac{w^n}{c_n} 
\sum_{k=0}^{\infty }\frac{z^k\overline w^k}{c_k^2}\,d\mu (w)\\
& = & \frac{1}{c_n}\, \int_D w^n 
\sum_{k=0}^{\infty   }\frac{z^k\overline w^{k+1}}{c_k^2}\,d\mu
(w)\\
& = & \frac{z^{n-1}}{c_n c_{n-1}^2}\int_D |w|^{2n}\,d\mu (w)\\
& = & \frac{c_n z^{n-1}}{c_{n-1}^2}, 
\end{eqnarray*}
where we used the fact that the series expansion of the Bergman
kernel converges uniformly on compact subsets (see \cite{K1} ).

Therefore we get
\begin{eqnarray*}
\|S(u_n)\|^2 & = &
\frac{1}{c_n^2}\int_D \left | \overline z\,z^n -
\frac{c_n^2}{c_{n-1}^2}\, z^{n-1} \right |^2 \,d\mu (z)\\
& = &
\frac{1}{c_n^2}\int_D      |z|^{2n-2}     \left     (     |z|^4
-\frac{2c_n^2|z|^2}{c_{n-1}^2}  + \frac{c_n^4}{c_{n-1}^4} \right
)\,d\mu (z)\\
& = &
\frac{1}{c_n^2}\int_D |z|^{2n+2}\,d\mu (z) -
\frac{2}{c_{n-1}^2}\int_D |z|^{2n}\,d\mu (z) +
\frac{c_n^2}{c_{n-1}^4}\int_D |z|^{2n-2}\,d\mu (z)\\
& = & \frac{c_{n+1}^2}{c_n^2} - \frac{c_n^2}{c_{n-1}^2}
\end{eqnarray*}

Hence
$$\sum_{n=0}^{\infty }\|S(u_n )\|^2 < \infty $$
if and only if
$$\lim_{n\to \infty}\frac{c_{n+1}^2}{c_n^2}< \infty .$$
\end{proof}

{\bf Example.}
First we consider the case where $D$ is the open unit disc
$\mathbb  D$  in  $\mathbb C$ and $d\mu (z) = (1-|z|^2)^{\alpha
}\,d\lambda (z), \ \alpha \ge 0.$ Here
$$c_n^2     =    \int_{\mathbb    D}|z|^{2n}\,(1-|z|^2)^{\alpha
}\,d\lambda (z) = \frac{\pi \, n!}{(\alpha +n+1)(\alpha +n)\dots 
(\alpha +1)},$$
hence
$$\frac{c_{n+1}^2}{c_n^2} = \frac{n+1}{\alpha + n + 2},$$
which     implies  by  Theorem 2.1 that  the  corresponding
canonical solution operator to $\overline \partial $ is always Hilbert Schmidt.
\vskip 0.5 cm

{\bf Example.}
We  now  consider  the  case  of the unit ball $\mathbb B^2$ in
$\mathbb  C^2$  and the weight 
$$d\mu (z) = (1-|z_1|^2-|z_2|^2)^{\alpha }
\,d\lambda (z) \
,\ \alpha \ge 0 .$$

Set
$$c_{n_1,n_2}^2      =     \int_{\mathbb     B^2}|z_1|^{2n_1}\,
|z_2|^{2n_2}\,(1-|z_1|^2-|z_2|^2)^{\alpha }\,d\lambda (z).$$

Using polar coordinates we get
$$c_{n_1,n_2}^2 = 4\pi ^2 \int_0^1 \, \int_0^{(1-r_2)^{1/2}}
r_1^{2n_1+1}\,      r_2^{2n_2+1}\,      (1-r_1^2-r_2^2)^{\alpha
}\,dr_1\,dr_2.$$

Now we substitute $s_1= 1-r_1^2-r_2^2\ , \ s_2=1-r_2^2$ and use
properties of the beta-function (see \cite{K1} )
\begin{eqnarray*}
c_{n_1,n_2}^2 & = &
\pi^2 \int_0^1\, \int_0^{s_2} (s_2-s_1)^{n_1}\, (1-s_2)^{n_2}\,
s_1^{\alpha }\, ds_1\,ds_2\\
& = & \pi^2 \, B(n_1+1, \alpha +1) \int_0^1 (1-s_2)^{n_2}\,
s_2^{\alpha + n_1 +1}\,ds_2\\
&  =  &  \pi^2  \, B(n_1+1, \alpha +1)\, B(n_2+1, \alpha +n_1 +
2)\\
&  =  &  \frac{\pi^2\,n_1!\,n_2!}{(\alpha  +  n_1+n_2+2)(\alpha
+n_1+n_2+1)\dots (\alpha +1)}
\end{eqnarray*}

The Bergman kernel has the form
\begin{eqnarray*}
K_{\alpha }(z,w) & = &
\frac{\alpha +1}{\pi^2 }\, \frac{1}{(1-z_1\overline w_1 -z_2
\overline w_2)^{\alpha +3}}\\
& = & \frac{1}{\pi^2} \, \sum_{n_1,n_2=0}^{\infty }
\frac{(\alpha +n_1+n_2+2)\dots (\alpha +1)}{n_1!\,n_2!}
z_1^{n_1}\,z_2^{n_2}\,\overline w_1^{n_1}\,
\overline w_2^{n_2}
\end{eqnarray*}

Now let 
$$u_{n_1,n_2}(z_1, z_2) =
\frac{z_1^{n_1}\,z_2^{n_2}}{c_{n_1,n_2}},$$
then the system 
$$\{ u_{n_1,n_2}(z_1, z_2)\,d\overline z_1 \ , \ 
u_{n_1,n_2}(z_1, z_2)\,d\overline z_2 \ : \ n_1, n_2\in \mathbb
N_0 \} $$
constitutes  an  orthonormal basis for the space $A_{(0,1)}^2(\mathbb B^2,
d\mu )$ of $(0,1)$-forms with coefficenits belonging to
$A^2(\mathbb B^2, d\mu ).$

Using the methods of Theorem 2.1 we obtain
\begin{eqnarray*}
\|S(u_{n_1,n_2}d\overline z_1) \|^2 & = &
\frac{c_{n_1+1,n_2}^2}{c_{n_1,n_2}^2}- 
\frac{c_{n_1,n_2}^2}{c_{n_1-1,n_2}^2}\\
& = &
\frac{\alpha + n_2+2}{(\alpha +n_1+n_2+3)(\alpha +n_1+n_2+2)}
\end{eqnarray*}
and
\begin{eqnarray*}
\|S(u_{n_1,n_2}d\overline z_2) \|^2 & = &
\frac{c_{n_1,n_2+1}^2}{c_{n_1,n_2}^2}- 
\frac{c_{n_1,n_2}^2}{c_{n_1,n_2-1}^2}\\
& = &
\frac{\alpha + n_1+2}{(\alpha +n_1+n_2+3)(\alpha +n_1+n_2+2)}.
\end{eqnarray*}

Since
$$\sum_{n_1, n_2 =1}^{\infty}[ 
\frac{\alpha +n_2+2}{(\alpha +n_1+n_2+2)(\alpha +n_1+n_2+3)}$$
$$+ \frac{\alpha +n_1+2}{(\alpha +n_1+n_2+2)(\alpha +n_1+n_2+3)} ]
=\infty ,$$
the canonical solution operator
$$S  :  A^2_{(0,1)}(\mathbb  B^2, d\mu )\longrightarrow L^2(\mathbb
B^2 ,d\mu )$$
fails to be a Hilbert-Schmidt operator. (For the case $\alpha =
0$ see \cite{Has2} ).

\section{Spaces of entire functions in one variable}

Now we consider weighted spaces on entire functions 
$$A^2(\mathbb C, e^{- |z|^m} )=\{ f: \mathbb C \longrightarrow \mathbb C \  
: \ \|f\|_m^2:= \int_{\mathbb C}|f(z)|^2 \, e^{-|z|^m} \,d\lambda (z)<\infty \},$$
where $ \ m>0.$
Let again 
$$c_k^2= \int_{\mathbb C}|z|^{2k} \, e^{-|z|^m} \,d\lambda (z).$$
Then 
$$K_m(z,w)= \sum_{k=0}^{\infty} \frac{z^k \overline w^k}{c_k^2}$$
is the reproducing kernel for 
$A^2(\mathbb C, e^{-|z|^m} ).$

In the sequel the expression 
$$\frac{c_{k+1}^2}{c_k^2} - \frac{c_k^2}{c_{k-1}^2} $$
will become important. Using the integral representation of the $\Gamma -$function 
one easily sees that the above expression is equal to
$$\frac{\Gamma \left ( \frac{2k+4}{m} \right ) }
{\Gamma \left ( \frac{2k+2}{m} \right )} -
\frac{\Gamma \left ( \frac{2k+2}{m} \right ) }{\Gamma \left ( \frac{2k}{m} \right )} .$$
For $m=2$ this expression equals to $1$ for each $k=1,2,\dots .$ 
We will be interested in the limit behavior for $k\to \infty .$ 
By Stirlings formula the limit behavior is equivalent
to the limit behavior of the expression
$$\left ( \frac{2k+2}{m} \right )^{2/m} - \left ( \frac{2k}{m} \right )^{2/m}, $$
as $k\to \infty .$ Hence we have shown the following 

\begin{lemm} The expression 
$$\frac{\Gamma \left ( \frac{2k+4}{m} \right ) }
{\Gamma \left ( \frac{2k+2}{m} \right )} -
\frac{\Gamma \left ( \frac{2k+2}{m} \right ) }{\Gamma \left ( \frac{2k}{m} \right )} $$
tends to $\infty $ for $0<m<2$, is equal to $1$ for $m=2$ and tends to zero for $m>2$ as 
$k$ tends to $\infty.$
\end{lemm}

Let $0 < \rho < 1,$ define $f_{\rho }(z):=f(\rho z )$ 
and $\tilde f_{\rho }(z) =  \overline z f_{\rho }(z),$
for $f\in A^2(\mathbb C, e^{-|z|^m}) .$
Then it is easily seen that $\tilde f_{\rho } \in L^2(\mathbb C, e^{-|z|^m}) ,$
but there are functions $g \in A^2(\mathbb C, e^{-|z|^m}) $ such that
$\overline z g \not\in  L^2(\mathbb C, e^{-|z|^m}) $.

Let $P_m :L^2(\mathbb C, e^{-|z|^m}) \longrightarrow A^2(\mathbb C, e^{-|z|^m}) $
denote the orthogonal projection. Then $P_m$ can be written in the form
$$P_m(f)(z)= \int_{\mathbb C }K_m(z,w)f(w)\,d\lambda (w) \ , \ f\in 
L^2(\mathbb C, e^{-|z|^m}) .$$

\begin{prop} 
Let $m\ge 2.$ Then there is a constant $C_m>0$ depending only on $m$ such that
$$\int_{\mathbb C}\left | \tilde f_{\rho }(z) - P_m(\tilde f_{\rho })(z) \right |^2
e^{-|z|^m}\,d\lambda (z) \le C_m \ \int_{\mathbb C}|f(z)|^2
e^{-|z|^m}\,d\lambda (z) ,$$
for each $0<\rho <1$ and for each $f\in A^2(\mathbb C, e^{-|z|^m}) .$
\end{prop}

\begin{proof} First we observe that for the Taylor expansion of
$f(z)=\sum_{k=0}^{\infty }a_kz^k $ we have
\begin{eqnarray*}
P_m(\tilde f_{\rho })(z) &=& \int_{\mathbb C}\sum_{k=0}^{\infty} 
\frac{z^k \overline w^k}{c_k^2}
\left (  \overline w \sum_{j=0}^{\infty }a_j \rho ^j w^j \right ) \ e^{-|w|^m }
\,d\lambda (w)\\ & = &
\sum_{k=1}^{\infty}a_k \, \frac{c_k^2}{c_{k-1}^2}\, \rho^k z^{k-1}.
\end{eqnarray*}
Now we obtain
\begin{eqnarray*}
\lefteqn{
\int_{\mathbb C}\left | \tilde f_{\rho }(z) - P_m(\tilde f_{\rho })(z) \right |^2
e^{-|z|^m}\,d\lambda (z) }\\ & = &
\int_{\mathbb C} \left ( \overline z \sum_{k=0}^{\infty}a_k \rho^k z^k -
\sum_{k=1}^{\infty}a_k \, \frac{c_k^2}{c_{k-1}^2}\, \rho^k z^{k-1}
\right )\\  & \times &
\left ( z \sum_{k=0}^{\infty}a_k \rho^k \overline z^k -
\sum_{k=1}^{\infty}\overline{a_k} \, \frac{c_k^2}{c_{k-1}^2}\, \rho^k
\overline z^{k-1}\right ) \ e^{-|z|^m }\,d\lambda (z)\\ &=&
\int_{\mathbb C } ( \sum_{k=0}^{\infty} |a_k|^2 \rho^{2k} |z|^{2k+2}
-2 \sum_{k=1}^{\infty }|a_k|^2\, \frac{c_k^2}{c_{k-1}^2}\, \rho^{2k}  |z|^{2k}\\ 
& + &
\sum_{k=1}^{\infty }|a_k|^2 \, \frac{c_k^4}{c_{k-1}^4}\,  \rho^{2k} |z|^{2k-2}
 )\ e^{-|z|^m }\,d\lambda (z)\\
&=& |a_0|^2\, c_1^2\,  +  \sum_{k=1}^{\infty }|a_k|^2\, 
c_k^2 \, \rho^{2k} \, \left ( \frac{c_{k+1}^2}{c_k^2} - \frac{c_k^2}{c_{k-1}^2} \right ).
\end{eqnarray*}

Now the result follows from the fact that
$$ \int_{\mathbb C}|f(z)|^2 e^{-|z|^m}\,d\lambda (z) =
\sum_{k=0}^{\infty }|a_k|^2\, c_k^2 ,$$
and that the sequence $\left ( \frac{c_{k+1}^2}{c_k^2} 
- \frac{c_k^2}{c_{k-1}^2} \right )_k$
is bounded. 
\end{proof}

By Fatou's theorem 
\begin{eqnarray*}
\lefteqn{
\int_{\mathbb C}\lim_{\rho \to 1 }
\left | \tilde f_{\rho }(z) - P_m(\tilde f_{\rho })(z) \right |^2
e^{-|z|^m}\,d\lambda (z)}\\ &\le &
\sup_{0<\rho <1}
\int_{\mathbb C}\left | \tilde f_{\rho }(z) - P_m(\tilde f_{\rho })(z) \right |^2
e^{-|z|^m}\,d\lambda (z)\\ & \le &
C_m \ \int_{\mathbb C}|f(z)|^2 e^{-|z|^m}\,d\lambda (z) 
\end{eqnarray*}
and hence the function 
$$F(z):= \overline z \sum_{k=0}^{\infty}a_k z^k -
\sum_{k=1}^{\infty}a_k \, \frac{c_k^2}{c_{k-1}^2}\, z^{k-1}$$
belongs to $L^2(\mathbb C, e^{-|z|^m}) $ and satiafies
$$\int_{\mathbb C} |F(z)|^2
e^{-|z|^m}\,d\lambda (z) \le C_m \ \int_{\mathbb C}|f(z)|^2
e^{-|z|^m}\,d\lambda (z) .$$
The above computation also shows that $\lim_{\rho \to 1}\|\tilde f_{\rho }-
P_m(\tilde f_{\rho })\|_m = \|F\|_m $ and by a standard argument for $L^p$-spaces
(see for instance \cite{E})
$$\lim_{\rho \to 1}\|\tilde f_{\rho }-P_m(\tilde f_{\rho })-F \|_m =0.$$

\begin{prop} 
Let $m\ge 2$ and consider an entire function $f\in A^2(\mathbb C, e^{-|z|^m}) $
with Taylor series expansion $f(z)=\sum_{k=0}^{\infty }a_kz^k .$ Let 
$$F(z):= \overline z \sum_{k=0}^{\infty}a_k z^k -
\sum_{k=1}^{\infty}a_k \, \frac{c_k^2}{c_{k-1}^2}\, z^{k-1}$$
and define $S_m(f):=F.$ Then $S_m:A^2(\mathbb C, e^{-|z|^m}) \longrightarrow 
L^2(\mathbb C, e^{-|z|^m}) $ is a continuous linear operator, representing the 
canonical solution operator to 
$\overline \partial $ restricted to $A^2(\mathbb C, e^{-|z|^m}) ,$
i.e. $\overline \partial S_m(f)=f$ and $S_m(f)\perp A^2(\mathbb C, e^{-|z|^m}) .$
\end{prop}

\begin{proof} 
A similar computation as in the proof of 
Proposition 1 in \cite{Has2}  shows that the function
$F$ defined above satisfies $\overline \partial F = f.$ 
Let $S_m(f):=F.$ Then, by the last remarks,
$S_m :A^2(\mathbb C, e^{-|z|^m}) \longrightarrow L^2(\mathbb C, e^{-|z|^m})$ 
is a continuous linear solution operator for $\overline \partial .$ 
For arbitrary $h\in A^2(\mathbb C, e^{-|z|^m})$ we
have
$$(h,S_m(f))_m=(h,F)_m=\lim_{\rho \to 1}(h,\tilde f_{\rho }-P_m(\tilde f_{\rho }))_m=
\lim_{\rho \to 1}(h-P_m(h),\tilde f_{\rho })_m=0,$$
where $(.\, ,.)_m$ denotes the inner product in $L^2(\mathbb C, e^{-|z|^m}) .$
Hence $S_m$ is the canonical solution operator for $\overline \partial $ restricted to
$A^2(\mathbb C, e^{-|z|^m}) .$
\end{proof}

\begin{theo} The canonical solution operator to $\overline \partial $ restricted to 
the space $A^2(\mathbb C, e^{-|z|^m}) $ is compact if and only if
$$\lim_{k\to \infty }\left ( \frac{c_{k+1}^2}{c_k^2} 
- \frac{c_k^2}{c_{k-1}^2} \right )=0.$$
\end{theo}

\begin{proof} 
For a complex polynomial $p$ the canonical solution operator $S_m$ can be written
in the form
$$S_m(p)(z)= \int_{\mathbb C}K_m(z,w)p(w)(\overline z - \overline w)\,d\lambda (w),$$
therefore we can express the conjugate $S_m^*$ in the form
$$S_m^*(q)(w)=\int_{\mathbb C}K_m(w,z)q(z)(z - w)\,d\lambda (z),$$
if $q$ is a finite linear combination of the terms $\overline z^k \, z^l.$
This follows by considering the inner product
$(S_m(p),q)_m=(p,S_m^*(q))_m .$

Now we claim that 
$$S_m^*S_m(u_n)(w)=\left ( \frac{c^2_{n+1}}{c^2_n}-\frac{c^2_n}{c^2_{n-1}}\right )
u_n(w)\ \ ,n=1,2,\dots $$
and
$$S_m^*S_m(u_0)(w)=\frac{c^2_1}{c^2_0}\ u_0(w),$$
where $\{ u_n(z)=z^n/c_n , k=0,1,\dots \}$ is the standard orthnormal basis of
$A^2(\mathbb C, e^{-|z|^m}).$ 

From the proof of Theorem 2.1 we know that
$$ S_m(u_n)(z)=\overline z u_n(z)-\frac{c_n z^{n-1}}{c_{n-1}^2}, \ n=1,2,\dots .$$
Hence 
\begin{eqnarray*}
S_m^*S_m(u_n)(w) &=&
\int_{\mathbb C} K_m(w,z) (z-w)\left ( \frac{\overline z z^n}{c_n}-
\frac{c_n z^{n-1}}{c_{n-1}^2}\right ) \,d\lambda (z)\\ &=&
\int_{\mathbb C} \sum_{k=0}^{\infty }\frac{w^k\overline z^k}{c_k^2}
 (z-w)\left ( \frac{\overline z z^n}{c_n}-
\frac{c_n z^{n-1}}{c_{n-1}^2}\right ) \,d\lambda (z). 
\end{eqnarray*}
This integral is computed in two steps: first the multiplication by $z$

$$
\int_{\mathbb C} \sum_{k=0}^{\infty }\frac{w^k\overline z^k}{c_k^2}
\left ( \frac{\overline z z^{n+1}}{c_n}-
\frac{c_n z^n}{c_{n-1}^2}\right ) \,d\lambda (z) $$
\begin{eqnarray*}
&=&
\int_{\mathbb C} \frac{z^{n+1}}{c_n}
\sum_{k=0}^{\infty }\frac{w^k\overline z^{k+1}}{c_k^2}\,d\lambda (z)-
\frac{c_n}{c^2_{n-1}}\int_{\mathbb C} z^n
\sum_{k=0}^{\infty }\frac{w^k\overline z^k}{c_k^2}\,d\lambda (z) \\ &=&
\frac{w^n}{c^3_n}\int_{\mathbb C} |z|^{2n+2}\,d\lambda (z)-
\frac{w^n}{c^2_{n-1}c^2_n}\int_{\mathbb C} |z|^{2n}\,d\lambda (z) \\ &=&
\left ( \frac{c^2_{n+1}}{c^3_n} - \frac{c_n}{c^2_{n-1}}\right )\ w^n
\end{eqnarray*}

And now the multiplication by $w$

$$w\int_{\mathbb C} \sum_{k=0}^{\infty }\frac{w^k\overline z^k}{c_k^2}
\left ( \frac{\overline z z^n}{c_n}-
\frac{c_n z^{n-1}}{c_{n-1}^2}\right ) \,d\lambda (z)$$
\begin{eqnarray*}
&=&
w\int_{\mathbb C} \frac{z^n}{c_n}
\sum_{k=0}^{\infty }\frac{w^k\overline z^{k+1}}{c_k^2}\,d\lambda (z)-
w\int_{\mathbb C} \frac{c_nz^{n-1}}{c^2_{n-1}}
\sum_{k=0}^{\infty }\frac{w^k\overline z^k}{c_k^2}\,d\lambda (z) \\ &=&
w\left ( \frac{c_n w^{n-1}}{c^2_{n-1}}-
\frac{c_n w^{n-1}}{c^2_{n-1}}\right ) \\ &=&
0,
\end{eqnarray*}
which implies that
$$S_m^*S_m(u_n)(w)=\left ( \frac{c^2_{n+1}}{c^2_n}-\frac{c^2_n}{c^2_{n-1}}\right )
u_n(w)\ \ ,n=1,2,\dots ,$$
the case $n=0$ follows from an analogous computation.

The last statement says that $S_m^*S_m$ is a diagonal operator with respect to the 
orthonormal basis $\{ u_n(z)=z^n/c_n\}$ of $A^2(\mathbb C, e^{-|z|^m}). $  
Therefore it is easily seen that $S_m^*S_m$ is compact if and only if
$$\lim_{n \to \infty }
\left ( \frac{c^2_{n+1}}{c^2_n}-\frac{c^2_n}{c^2_{n-1}}\right )=0.$$
Now the conclusion follows, since $S_m^*S_m$ is compact if and only if $S$ is compact
(see for instance \cite{Wei}). 
\end{proof}

\begin{theo} The canonical solution operator for $\overline \partial $ restricted to
the space $A^2(\mathbb C, e^{-|z|^m}) $ is compact, if $m>2.$ 
The canonical solution operator
for $\overline \partial $ as operator from 
$L^2(\mathbb C, e^{-|z|^2}) $ into itself is not compact.
\end{theo}

\begin{proof} The first statement follows immediately from Theorem 3.1 and Lemma 3.1. 
For the
second statement we use H\"ormander's $L^2 $-estimate 
for the solution of the $\overline \partial $
equation \cite{H} : for each $f\in L^2(\mathbb C, e^{-|z|^2}) $ there is a function
$u\in L^2(\mathbb C, e^{-|z|^2}) $ such that $\overline \partial u =f $ and
$$\int_{\mathbb C}|u(z)|^2 \, e^{-|z|^2} \,d\lambda (z) \le 4 
\int_{\mathbb C}|f(z)|^2 \, e^{-|z|^2} \,d\lambda (z).$$
Hence the canonical solution operator for $\overline \partial $ as operator from
$L^2(\mathbb C, e^{-|z|^2}) $ into itself is continuous and its restriction
to the closed subspace $A^2(\mathbb C, e^{-|z|^2}) $ fails to be compact by
Propositon 3.2 and Lemma 3.1. By the definition of compactness this implies that the
canonical solution operator is not compact 
as operator from $L^2(\mathbb C, e^{-|z|^2}) $ 
into itself. 
\end{proof}

{\bf Remark.}
In the case of the Fock space $A^2(\mathbb C, e^{-|z|^2}) $ the compostion
$S_2^*S_2 $ equals to the identity on $A^2(\mathbb C, e^{-|z|^2}), $ which follows 
from the proof of Theorem 3.1.

\begin{theo} Let $m\ge 2.$ The canonical solution 
operator for $\overline \partial $ restricted to
$A^2(\mathbb C, e^{-|z|^m}) $ fails to be Hilbert Schmidt.
\end{theo}

\begin{proof} By Proposition 3.2 we know 
that the canonical solution operator is continuous
and we can apply the method from Theorem 2.1 : now in our case we have
$$\frac{c_{n+1}^2}{c_n^2}  = \Gamma \left
( \frac{2n+4}{m} \right ) / \Gamma \left
( \frac{2n+2}{m} \right ),$$
which,  by  Theorem 2.1 and Stirling's  formula, implies that the corresponding
canonical  solution  operator  to $\overline \partial $ fails to be Hilbert
Schmidt. 
\end{proof}

In the case of several variables the corresponding operator $S^*S$ is more complicated,
nevertheless we can handle a sligthly different situation with different
methods from functional analysis (see next section).

\section{Weighted spaces of entire functions in several variables}

In  this  part  we show that the canonical solution operator to
$\overline \partial  $  is  a  Hilbert-Schmidt operator for a wide class of
weighted spaces of entire functions.

The weight functions we are considering are of the form $z \mapsto
\tau p(z),$ where $\tau > 0$ and $p: \mathbb C^n \longrightarrow
\mathbb R.$ We suppose that $p $ is a plurisubharmonic function satisfying
$$p^*(w) := \sup \{ \Re <z,w> - p(z) \ : \ z\in \mathbb C^n \}
< \infty .$$

Then $p^{**}=p$ and
$$\lim_{|z|\to \infty}\frac{p(z)}{|z|} = \infty $$
(see Lemma 1.1. in \cite{Has1}). And it is easily seen that 
$$\int_{\mathbb C^n }\exp [(\tau - \sigma )\,p(z)]\,d\lambda (z)
< \infty ,$$
whenever $\tau - \sigma < 0.$

We further assume that 
$$\lim_{|z|\to \infty }\frac{\tilde p (z)}{p(z)} = 1,$$
where $\tilde p (z) = \sup \{ p(z+\zeta ) \, : \, |\zeta |\le 1 \}.$

It follows that the last property is equivalent to the following
condition: for each $\sigma >0$ and for each $\tau >0 $ with
$\tau < \sigma $ there is a constant $C=C(\sigma, \tau )>0$ such 
that
$$\tau \, \tilde p(z) - \sigma \, p(z) \le C ,$$
for each $z\in \mathbb C^n .$

Let $A^2 (\mathbb C^n, \sigma p)$ denote the Hilbert space of all entire 
functions $h: \mathbb C^n \longrightarrow \mathbb C $ such that
$$\int_{\mathbb C^n } |h(z)|^2 \, \exp (-2\sigma p(z))\,
d\lambda (z) < \infty .$$

\begin{theo} Suppose that $p$ is a weight function with the
properties listed above. Then for each $\sigma >0 $ there exists 
a number $\tau > 0 $ with $\tau < \sigma $ such that the canonical
solution operator $S_1 $ to $\overline \partial $ is a Hilbert-Schmidt operator
as a mapping
$$S_1 : A^2_{(0,1)} (\mathbb C^n, \tau p) \longrightarrow
L^2(\mathbb C^n, \sigma p).$$
\end{theo}

\begin{proof} By Lemma 28.2 from \cite{MV} we have to show that 
$$\left [
\int_{\mathbb C^n}|S_1(f)(z)|^2\, \exp (-2\sigma p(z))\,d\lambda (z)
\right ]^{1/2}
\le \int_{\mathcal U}|(f,g)|\, d\mu (g), $$
where $(.,.)$ denotes the inner product of the Hilbert space
$A^2_{(0,1)} (\mathbb C^n, \tau p),$ \ $ \mathcal U$ is the unit ball of
$A^2_{(0,1)} (\mathbb C^n, \tau p)$ \ , \  $\mu $ is a  Radon
measure   on   the   weakly   compact   set  $\mathcal  U$  and
$f=\sum_{j=1}^n f_j\,d\overline z_j $ and
$g=\sum_{j=1}^n g_j\,d\overline z_j .$
\vskip 0.5 cm
We  first  show  that for $0< \tau < \tau_1 < \tau_2 < \tau_3 <
\sigma $ we have
\begin{eqnarray*}
\lefteqn{
\left [ \int_{\mathbb  C^n}|f_j(z)|^2\,  \exp  (-2\tau_3  p(z))\,d\lambda
(z) \right ]^{1/2} }\\
& \le &
C_{\tau_3 ,\tau_2 }\, \sup \{ |f_j(z)|\, \exp (-\tau_2 p(z)) \, : \, z\in 
\mathbb C^n \}\\
& \le &
C_{\tau_2 , \tau_1 }\, 
\int_{\mathbb C^n}|f_j(z)|\, \exp (-\tau_1 p(z))\,d\lambda (z)\\
& \le &
C_{\tau_1 , \tau }
\left [ \int_{\mathbb C^n}|f_j(z)|^2\, \exp (-2\tau p(z))\,d\lambda (z) 
\right ]^{1/2},\\
\end{eqnarray*}
for each $f \in A^2_{(0,1)} (\mathbb C^n, \tau p).$

To  show  this assertion we make use of the assumption that the
coefficients of the $(0,1)$-form $f$ are entire functions:

The first inequality follows from the fact that
$$\int_{\mathbb  C^n} \exp ((2\tau_2 - 2\tau_3 )p(z))\,d\lambda (z) <
\infty .$$

For   the  second  inequality  use  Cauchy's  theorem  for  the
coefficients  $f_j$ of $f$ to show that for $B_z = \{ \zeta \in
\mathbb C^n \, : \, |z-\zeta |\le 1 \}$ we have
\begin{eqnarray*}
\lefteqn{
|f_j(z)|  \le  C  \, \int_{B_z}|f_j(\zeta )|\,d\lambda (\zeta )}\\
& = & 
C\, \int_{B_z}|f_j(\zeta )|\exp ( -\tau_1 p(\zeta ))
 \exp (\tau_1 p(\zeta ))\,d\lambda (\zeta )\\
 & \le &
C\, \int_{B_z}|f_j(\zeta )| \exp ( -\tau_1 p(\zeta ))
 \,d\lambda (\zeta ) \ 
\sup \{ \exp (\tau_1 p(\zeta )) \, : \, \zeta \in B_z \}\\
& \le &
C'\, \int_{\mathbb C^n}|f_j(\zeta )| \exp ( -\tau_1 p(\zeta ))
 \,d\lambda (\zeta ) \ \exp (\tau_2 p(z)),\\
 \end{eqnarray*}
where we used the properties of the weight function $p.$

The  third  inequality  is  a consequence of the Cauchy-Schwarz
inequality:
\begin{eqnarray*}
\lefteqn{
\int_{\mathbb C^n}|f_j(\zeta )| \exp ( -\tau_1 p(\zeta ))
 \,d\lambda (\zeta ) }\\
 & = &
\int_{\mathbb C^n}|f_j(\zeta )| \exp ( -\tau p(\zeta ))
\exp ((\tau - \tau_1)p(\zeta ))
 \,d\lambda (\zeta ) \\
 & \le &
 \left [ \int_{\mathbb C^n}|f_j(\zeta )|^2 \exp ( -2\tau p(\zeta ))
 \,d\lambda (\zeta ) \right ]^{1/2} \,
 \left [ \int_{\mathbb C^n} \exp ( (2\tau - 2\tau_1 )p(\zeta ))
 \,d\lambda (\zeta ) \right ]^{1/2}.\\ 
 \end{eqnarray*}
By  H\"ormander's  $L^2$-estimates (\cite{H}, Theorem 4.4.2]) we have
for $\tau  <  \tau_3  <  \sigma  $ and the properties of the weight
function
\begin{eqnarray*}
\lefteqn{
\int_{\mathbb  C^n}|S_1(f)(z)|^2\,  \exp (-2\sigma p(z))\,d\lambda
(z)}\\
& \le &
\int_{\mathbb C^n}|S_1(f)(z)|^2\, \exp (-2\tau_3 p(z))
\, (1+|z|^2)^{-2}\,d\lambda (z)\\
& \le &
\int_{\mathbb  C^n}|f(z)|^2\,  \exp  (-2\tau_3  p(z))\,d\lambda
(z).\\
\end{eqnarray*}
Here  we  used  the  fact  that the canonical solution operator
$S_1$ can be written in the form $S_1(f)=v-P(v) ,$ where $v$ is 
an  arbitrary  solution  to  $\overline \partial  u  =  f$ belonging to the
corresponding Hilbert space and that $\|S_1(f)\|=\|v-P(v)\|=
\min \{ \|v-h\| \, : \, h\in A^2 \}\le \|v \|.$

Now choose $\tau_2 $ such that 
$\tau    <  \tau_2  < \tau_3  <  \sigma  ,$ then we obtain from
the above inequalities
\begin{eqnarray*}
\lefteqn{
\left [
\int_{\mathbb  C^n}|S_1(f)(z)|^2\,  \exp (-2\sigma p(z))\,d\lambda
(z) \right ]^{1/2}}\\
& \le &
D\, \int_{\mathbb  C^n}|f(z)|_1\,  \exp  (-\tau_2  p(z))\,d\lambda
(z),\\
\end{eqnarray*}
where  $D>0$  is  a  constant and $|f(z)|_1 := |f_1(z)|+\dots +
|f_n(z)|.$

Now define for $z\in \mathbb C^n$ and $\tau < \tau_1 < \tau_2 $
$$\delta_z^{\tau_1, \tau }(f_j) :=  C_{\tau_1 , \tau }^{-1}\,
f_j(z)\, \exp (-\tau_1 p(z)).$$                          
Then
\begin{eqnarray*}
\lefteqn{
 C_{\tau_1 , \tau }^{-1}\, \sup \{ |f_j(z)|\, \exp (-\tau_1 p(z))
\, : \, z\in \mathbb C^n \} }\\
& \le &
\left [ \int_{\mathbb C^n } |f_j(z)|^2\, \exp (-2\tau p(z))\,
d\lambda (z) \right ]^{1/2},\\
\end{eqnarray*}
which, by the Riesz representation theorem for the Hilbert
space $A^2_{(0,1)} (\mathbb C^n, \tau p )$, means that each 
$\delta_z^{\tau_1,  \tau  }$  can  be  viewed  as  an element of
$\mathcal U .$

For $\phi \in \mathcal C (\mathcal U )$ the expression
$$\mu  (\phi ) = D \,\int_{\mathbb C^n} \phi (\delta_z^{\tau_1,
\tau })\,\exp ((\tau_1 - \tau_2 )p(z))\,d\lambda (z)$$
defines a Radon measure on the weakly compact set $\mathcal U.$

This follows from the fact that
$$\mu  (\phi ) \le D \, \sup \{ |\phi (g)| \, :\, g\in \mathcal
U \}\ 
\int_{\mathbb C^n} \exp ((\tau_1 - \tau_2 )p(z))\,d\lambda (z).$$

Now  take  for  $\phi $ the continous functions $\phi_j (g_j) =
|(f_j,g_j)|,$ where $f_j $ is fixed. Then
\begin{eqnarray*}
\lefteqn{
\int_{\mathcal U}|(f,g)|\, d\mu (g)}\\
& = &
 D \,\int_{\mathbb C^n} |f(z)|_1\, \exp (-\tau_1 p(z))
 \,\exp ((\tau_1 - \tau_2 )p(z))\,d\lambda (z)\\
 & = &
D \,\int_{\mathbb C^n} |f(z)|_1\, 
\exp (- \tau_2 p(z))\,d\lambda (z)\\
\end{eqnarray*}
and hence
$$\left [
\int_{\mathbb C^n}|S_1(f)(z)|^2\, \exp (-2\sigma p(z))\,d\lambda (z)
\right ]^{1/2}
\le  \int_{\mathcal  U}|(f,g)|\, d\mu (g). $$
\end{proof}

\end{document}